\documentclass[12pt]{amsart}
\usepackage{epsfds} 

\usepackage{epic, eepic}
\def\squarebox#1{\hbox to #1{\hfill\vbox to #1{\vfill}}} 
\newcommand{\stopthm}{\hfill\hfill\vbox{\hrule\hbox{\vrule\squarebox 
                 {.667em}\vrule}\hrule}\smallskip}

\newcommand{\cP}{{\mathcal P}} 
\newcommand{\CC}{{\mathbb C}}
\newcommand{\pa}{\partial}
\newcommand{\CI}{{\mathcal C}^\infty }
\newcommand{\DI}{{\dot {\mathcal C}}^\infty }
\newcommand{\Oo}{{\mathcal O}} 

\newcommand{\E}{{\mathcal E}}
\newcommand{\up}{\Upsilon}
\newcommand{\ep}{\epsilon}

\newcommand{\G}{{\Gamma}}
\newcommand{\cG}{{\mathcal G}}

\newcommand{\RR}{{\mathbb R}}
\newcommand{\NN}{{\mathbb N}}
\newcommand{\HH}{{\mathbb H}}

\newcommand{\N}{{\mathbb N}}

\newcommand{\vol}{\operatorname{vol}}

\newcommand{\Ric}{\operatorname{Ric}}
\newcommand{\Res}{\operatorname{Res}}

\newcommand{\Diff}{\operatorname{Diff}}
\newcommand{\tr}{\operatorname{tr}}
\newcommand{\Xint}{\buildrel \circ\over X}

\renewcommand{\Re}{\mathop{\rm Re}\nolimits}

\theoremstyle{plain}

\newtheorem{thm}{Theorem}
\newtheorem{prop}{Proposition}[section]

\theoremstyle{definition}

\newtheorem*{rem}{Remark}

\numberwithin{equation}{section}

\title{Scattering matrix in conformal geometry}
\author{C. Robin Graham}
\address{Department of Mathematics, University of Washington,
Box 354350\\
Seattle, WA 98195}
\email{robin@math.washington.edu}
\author{Maciej Zworski}
\address{Department of Mathematics, University of California\\
Berkeley, CA 94720}
\email{zworski@math.berkeley.edu}

\begin{document}
\maketitle
\section{Statement of the results}

This paper describes the connection between
scattering matrices on conformally compact asymptotically Einstein
manifolds and conformally invariant objects on their boundaries at
infinity.  This connection is a manifestation of the general principle
that the far field phenomena on a conformally compact 
Einstein manifold are related to conformal theories on its boundary at
infinity.  This relationship was proposed in 
\cite{FG} as a means of studying conformal geometry, and the principle
forms the basis of 
the AdS/CFT correspondence in quantum gravity -- see
\cite{Mal},\cite{W},\cite{HS},\cite{GW} and references given there. 

We first define the basic objects discussed here.
By a conformal structure on a compact manifold $ M$ we mean 
an equivalence class $ [h] $ determined by a metric representative $h$: 
\[ \hat{h} \in [h] \ \Longleftrightarrow \ \hat{h} = e^{ 2 \Upsilon } h \,, \ \ 
\Upsilon \in \CI ( M ) \,.\]
Let $X$ be a compact $n+1$-manifold
with $ \partial X = M $, or in the case $ M$ does not bound a manifold, with
$ \partial X = M \sqcup M $. This can be achieved trivially by putting
$ X = [ 0 , 1 ] \times M $. 

Let $ x $ be a defining function of 
$ \partial X  $ in $ X$:
\[ x|_{\Xint}>0 \,, \ \ x |_{\partial X } = 0 \,, \ \ 
dx |_{ \partial X } \neq 0 \,. \]
We say that $ g $ is a 
conformally compact metric on $ X $ with conformal infinity
$ [h] $ if
\begin{equation}
\label{eq:confc}
 g = \frac {\overline g }{ x^2 } \,, \ \ {\overline g}|_{T\partial X} \in
[h] \,,\end{equation}
where ${\overline g}$ is a smooth metric on $X$. Since we can choose
different defining functions, the metric $ g$ determines only the conformal 
class $ [h]$. A conformally compact metric is said to be asymptotically
hyperbolic if its sectional curvatures approach $-1$ at $\pa X$; this is
equivalent to $|dx|_{\bar g}=1$ on $\pa X$.  
The basic example is hyperbolic space $
\HH^{n+1} $, with boundary given by $ \RR^{n} $ in the half-space
model and by $ {\mathbb S}^n $ in the ball model: the two being
conformally equivalent.  

One of the results of \cite{FG} is that given a conformal structure 
$[h]$ on $M$, one can construct a conformally
compact metric $g$ with conformal infinity $[h]$ which 
satisfies
\begin{equation}\label{ae}
 {\text{Ric}} (g ) + ng = \left\{ 
\begin{array}{ll} \Oo ( x^\infty  ) & \text{ for $ n $ odd} 
\\ \Oo ( x^{{n-2} } ) &  \text{ for $ n $ even.} 
\end{array} \right.
\end{equation}
When $n$ is even, the condition (\ref{ae}) is augmented by a
vanishing trace condition to the next order.  We call a metric $g$
satisfying these conditions asymptotically Einstein.
When $n$ is odd, the condition (\ref{ae}) 
together with an asymptotic evenness condition uniquely determine $g$ mod 
$ \Oo (x^{\infty})$ up to diffeomorphism.  We shall call a
metric $g$ 
which is asymptotically Einstein and which also satisfies 
the asymptotic evenness condition, a Poincar\'e metric associated to
$[h]$. 

Our first theorem relates the scattering matrix of a Poincar\'e
metric $g$ to the ``conformally
invariant powers of the Laplacian'' on $(M,[h])$.  The 
scattering matrix of $(X,g)$ is a meromorphic family $S(s)$ of
pseudodifferential operators on $M$ defined in terms of the
behaviour at infinity of solutions of $[\Delta_g -s(n-s)]u=0$,
which we discuss in \S \ref{scat}.  
The conformally invariant powers of the Laplacian are a family $P_k$, 
$k\in \NN$ and $k\leq n/2$ if $n$ is even, of 
scalar differential operators on $M$ constructed in \cite{GJMS}, which we
discuss in \S \ref{conf}.
These operators are {\em natural} in the sense that they can be written in
terms of covariant derivatives and curvature of a representative metric
$h$, and they are {\em invariant}
in the sense that if $\hat{h} = e^{ 2 \Upsilon } h$, then 
\begin{equation}
\label{eq:1.inv}
\widehat{P} _k =  e^{ ( - n/2 - k ) \Upsilon } 
P_k\, e^{ ( n/2 - k ) \Upsilon } \,. 
\end{equation}
The operator $P_k$ has the same principal part as $\Delta^k$ (in our
convention the Laplacian is a positive operator) and equals $\Delta^k$ if
$h$ is flat.

\begin{thm}
\label{t:1}
Let $ (M^n , [h] ) $ be a compact manifold with a conformal structure, 
and let $ ( X, g ) $ be a Poincar\'e
metric associated to $[h]$. Suppose that  
$  k \in \NN$ and $k\leq n/2$ if $n$ is even,
and that $ (n/2)^2 - k^2 $ is {\em not} an $L^2$-eigenvalue of $ \Delta_g
$.  If $ S ( s )$ is the scattering matrix of $ ( X, g ) $, and $ P_k $
the  conformally invariant operator on $M$, then $S(s)$ has a simple
pole at $s=n/2+k$ and 
\begin{equation}
\label{eq:t.1}
c_k  P_k  = - {\rm{Res}\;}_{ s = n/2 + k } \; S ( s) \,, 
\ \ c_{ k } = (-1)^{k} [ 2^{2k} k! ( k-1) ! ]^{-1} \,,
\end{equation}
where $ {\rm{Res}\;}_{ s = s_0 } \; S ( s) $ denotes the
residue at $s_0$ of the meromorphic family of operators $ S ( s )$.
\end{thm}

We remark that the condition on the spectrum
is automatically satisfied if $k\geq n/2$ and in general can be guaranteed
by perturbing the metric $ g $ in the interior.
Since $S(s)$ is self-adjoint for $s\in \RR$, as an immediate consequence
of Theorem~\ref{t:1} and of this remark we obtain 

\medskip

\noindent
{\bf Corollary.}
{\em The conformally invariant operators, $ P_k $, are self-adjoint.}

\medskip

This was previously known only for small values of $k$ for which the
operators can be explicitly calculated.
The first two of the invariant operators are the conformal Laplacian
\[
P_1 = \Delta + \frac{n-2}{4(n-1)} R 
\]
and the Paneitz operator
\[
P_2
= \Delta^2 + \delta Td +(n-4)(\Delta J +\frac{n}{2}J^2 -2|P|^2)/2. 
\]
Here $R$ denotes the scalar curvature, $J=R/(2(n-1))$, 
$P_{ij}=\frac{1}{n-2}(R_{ij}-Jh_{ij})$ where $R_{ij}$ is the Ricci
curvature, $T=(n-2)Jh -4P$ acting as an endomorphism on 1-forms, 
$|P|^2=P_{ij}P^{ij}$, and $\delta$ is the adjoint of $d$ (the
divergence operator). 

Another important notion of conformal geometry 
is Branson's $ Q $-curvature in even dimensions.  It
is a scalar function on $M$ constructed from the curvature tensor and its 
covariant derivatives, 
with an invariance property generalizing that of scalar curvature 
in dimension two:  if once again $\hat{h} = e^{ 2 \Upsilon } h$, then 
\begin{equation}
\label{eq:1.q}
e^{ n \Upsilon } {\widehat Q} = Q + P_{n/2}\Upsilon \,.
\end{equation}
There has been great progress recently in understanding 
the $Q$-curvature and its geometric meaning in low dimensions and on
conformally flat manifolds -- see
\cite{AC} for an example of recent work.  However, in general it remains a
rather mysterious quantity -- its definition (given in \cite{B} and
reviewed in \S \ref{conf} below)
in the general case is via analytic continuation in the dimension.  
In dimension two it is given by $Q=R/2$, and 
in dimension four by
$6 Q =  \Delta R + R^2 - 3 | {\rm{Ric}}|^2$.

If $n$ is even, then the operator $P_{n/2}$ has no constant term, i.e.
$P_{n/2}1=0$.  It therefore follows from Theorem~\ref{t:1} that $S(s)1$
extends holomorphically across $s=n$, so $S(n)1=\lim_{s\rightarrow
n}S(s)1$ is a well-defined function on $M$.
\begin{thm}
\label{t:2}
With the notation of Theorem \ref{t:1}, for $ n$ even, we have 
\begin{equation}
\label{eq:t.2}
c_{n/2} Q = S(n ) 1 \,. 
\end{equation}
\end{thm}

Theorem \ref{t:2} can be used as an alternative definition of the 
$ Q$-curvature.  We show in \S 4 that the conformal transformation law
(\ref{eq:1.q}) is an easy consequence of Theorems \ref{t:1} and \ref{t:2}.   
It follows from (\ref{eq:1.q}), the self-adjointness of $P_{n/2}$, and the
fact that $P_{n/2}1=0$, that $ \int_M Q $ is a conformal invariant.
For $(M,[h])$ conformally flat, it follows from a result of \cite{BGP} that
$\int_M Q$ is a multiple of the Euler characteristic $\chi(M)$.

A specific mathematical object which appeared in the study of the AdS/CFT 
correspondence is the renormalized volume of an asymptotically hyperbolic 
manifold $(X,g)$ -- see \cite{G} for a discussion and references.
It has also appeared in geometric scattering theory -- see
\cite{GZ},\cite{JS2}. 
As shown in \cite{GL} (Lemma 5.2 and the subsequent paragraph), a
choice of metric $h$ in the conformal class on $\pa X$
uniquely determines a defining function $x$ near $\pa X$ and an
identification of a neighborhood of $\pa X$ with $\pa X \times
[0,\epsilon]$ such that $g$ takes the form 
\begin{equation}\label{normalform}
g=x^{-2}(h_x + dx^2),\qquad h_0=h,
\end{equation}
where $h_x$ is a 1-parameter family of metrics on $\pa X$.
The renormalized volume is 
defined as the finite part in the expansion of
$ \vol_g ( \{ x > \epsilon \} ) $ as $ \epsilon \rightarrow 0 $.
For asymptotically Einstein metrics the expansions take 
a special form 
\begin{equation}
\label{eq:1.L}
 \begin{split}
&  \vol_g ( \{ x > \epsilon \} ) = c_0 \epsilon^{-n} + c_2 \epsilon^{-n+2}
+ \cdots + c_{n-1} \epsilon^{-1} + V + o(1) \\ 
& \ \\
& \ \ \ \ \ \ \text{for  $ n$  odd, }\\
& \ \\
&  \vol_g ( \{ x > \epsilon \} ) = c_0 \epsilon^{-n} + c_2 \epsilon^{-n+2}
+ \cdots + c_{n-2} \epsilon^{-2} 
 +  L \log ( 1 / \epsilon ) 
+ V + o(1) \\  
& \ \\
& \ \ \ \ \ \ \text{for $ n$ even. }
\end{split}
\end{equation}
 
It turns out that for asymptotically Einstein metrics $g$,
$V$ is independent of
the conformal representative $h$ on the boundary at infinity when $n$ is
odd, and $ L$ is independent of the conformal representative when $n$ is
even.  The dependence of $ V $ on the choice of $h$ for $n$ even is the
so-called holographic anomaly -- see \cite{HS},\cite{G}.  
Anderson \cite{A} has recently identified $V$ when $n=3$.
In an appendix to \cite{PP}, Epstein shows that for conformally compact
hyperbolic manifolds,
the invariants $L$ for $n$ even and $V$ for $n$ odd are each
multiples of the Euler characteristic $\chi(X)$.

Using the connection with the scattering matrix, we are 
able to identify $ L$ in terms of the $ Q$-curvature:
\begin{thm}
\label{t:3}
Let $ n $ be even and let $ L$ be defined by \eqref{eq:1.L}. Then
\begin{equation}
\label{eq:t.3}
L = 2 c_{n/2} \int_M Q \,,
\end{equation}
where $ c_{n/2} $ is defined in \eqref{eq:t.1}.
\end{thm}

We should stress that despite the fact that the scattering matrix is a
global object, in some sense
our results are all formal Taylor series statements at the boundary of 
$X$.  In fact, in \cite{FG2} it will be shown that, using a variant of the
ideas introduced here, a direct definition of $Q$ and proofs of 
Theorem \ref{t:3} and the self-adjointness of the $ P_k$'s
can be given based purely on formal asymptotics, 
avoiding the analytic continuation via the scattering matrix.
It is nevertheless worthwhile to proceed with the full scattering theory:
as a byproduct, this allows us to clarify certain 
confusing issues about the infinite rank 
poles of the scattering matrix at $ s = n/2 + l/2 $, $ l \in \NN$.

The paper is organized as follows. In \S \ref{ex} we present a 
simple one dimensional 
introduction to the relevant aspects of scattering theory. 
The more involved theory for asymptotically hyperbolic manifolds 
is then discussed in \S \ref{scat}: using results of Mazzeo-Melrose
\cite{MM} we give direct arguments for 
the existence and properties of the Poisson operator and the scattering
matrix, focussing particularly on their behaviour for $s$ near $n/2 +
\NN/2$. Although it is not relevant for our main results, 
at the end of \S \ref{scat} we outline 
a method for the study of the structure of the scattering matrix 
complementing the treatment in \cite{JS}. 
In \S \ref{conf} we show how the invariant operators $P_k$ may be
constructed from a Poincar\'e metric and discuss Branson's $Q$-curvature.  
Finally, \S \ref{pr} combines scattering theory with \S \ref{conf},
providing proofs of the main results.

In the paper, $ \NN = \{ 1, 2, \cdots \} $ denotes the natural numbers and  
$ \NN_0 = \NN \cup \{ 0 \} $. Linear operators are identified with their
distributional (integral) kernels.

\medskip
\noindent
{\sc Acknowledgments.} We would like to thank Rafe Mazzeo and Andras
Vasy for helpful comments.  The first author thanks the Mathematical 
Sciences Research Institute, where this project was started.  The second 
author thanks the Erwin Schr\"odinger Institute, where it was continued,
and the National Science Foundation for 
partial support under the grant DMS-9970614.

\section{A simple example}
\label{ex}

As a simple-minded introduction to scattering and as an illustration 
of the connection between residues of the scattering matrix and asymptotic
expansions we present a one-dimensional example. 

Thus we consider scattering on a
half line, $ X = [0, \infty ) $ with a compactly supported 
real valued potential
$ V \in L^\infty _{\rm{comp}} ( X; \RR )$. The quantum Hamiltonian is 
given by $ H = - \partial_y^2 + V( y ) $ and we impose, say, the 
Dirichlet boundary condition at $ y = 0 $. We are interested in 
the properties of generalized eigenfunctions at energies $ - s^2 $:
\[
( H + s^2 ) u ( y ) = 0 \,, \ \ u ( 0 ) = 0 \,, 
\]
which for large values of $ y $ and for $s\neq 0$ satisfy
\begin{equation}
\label{eq:2.1}
u ( y ) = A (s ) e^{s y  } + B (s) e^{-s y }.
\end{equation}
For $s\in i \RR \setminus \{0\}$, consideration of the Wronskian of
$ u $ and $ \bar u $ shows that $|A( s) |^ 2 = |B( s) |^2$.
Normalizing, we define the scattering matrix,
$ S ( s )$ (which in this case is a one-by-one matrix!) by 
\[
 S ( s ) = \frac{ B(s ) }{ A ( s) } \,. \]
This defines a meromorphic function of $ s \in \CC$, regular 
for $s\in i\RR$. 
The definition shows that
$ S ( s) = S ( -s )^{-1} $ and that for $ s \in i \RR$,  
$ S ( s )^* = S ( s)^{-1} $. Hence for all $ s $ we have 
\begin{equation}
\label{eq:2.unit} 
 S (- \bar s )^*  = S ( s )^{-1} \,, \  \  S ( \bar s )^* = S ( s ) \,,
\end{equation}
so that, in particular, $ S ( s ) $ is self-adjoint for $ s $ real.
When $ V $ is compactly supported then $L^2$-eigenvalues of $ H $, 
$ E \leq 0 $, correspond to the poles of $ S ( s )$, $ E = - s^2 $, 
$ \Re s > 0 $. This follows easily from \eqref{eq:2.1} and 
self-adjointness of $ H $. The poles of $ S $ for $ \Re s < 0 $ 
correspond to {\em resonances} -- see \cite{ZwN}.

When $ V $ is not compactly supported, then under relatively mild
assumptions we still have the scattering matrix for 
$ s \in i \RR $ but
its meromorphic continuation becomes very sensitive to the behaviour 
of $ V $ at infinity. The first physical case in which these 
difficulties\footnote{This is potentially rather confusing. The ``false''
poles of the scattering matrix discussed below almost led Heisenberg to 
abandoning his $ S$-matrix formalism, until a clear explanation was
provided by Jost -- see \cite{N}.}
occurred was that of Yukawa potential, $  V( y ) = e^{-y } $. 

We can study the Yukawa potential scattering using simple regular singular
point analysis, not unlike what one encounters in the study of the free
hyperbolic space. We start by making a change of variables:
\[ x = e^{-y}  \,, \ \ H + s^2 = -  ( x \partial_x)^2 + x + s^2 \,, \ \ 
X = (0, 1] \,. \]
Note that the definition of the scattering matrix in the new variables is
\begin{gather}
\label{eq:2.scat}
\begin{gathered}
( H + s^2 ) u ( x ) = 0 \,, \ \ u ( 1 ) = 0 \,, 
\ \ s \in i \RR \setminus \{0\}\,, \\
 u ( x ) = x^{-s } + S ( s) x^s + \Oo ( x )\,, \ \ x \rightarrow 0 \,.
\end{gathered}
\end{gather}
We obtain $ u (x) $ from the 
two solutions of $ ( H + s^2 ) G = 0 $:
\begin{equation}
\label{eq:2.3}
G_{ \pm } ( x , s ) = x^{ \pm s } \sum_{ j = 0 }^\infty b_j^\pm ( s) x^j \,,
\ \ b_j ^{\pm } ( s) = [ j! \Gamma ( \pm 2 s + j + 1 ) ]^{-1} \,. 
\end{equation}
These are independent provided $ 2 s \notin {\mathbb Z} $. 
The scattering matrix is obtained by finding a combination matching the
boundary conditions:
\[ S  ( s) = - \frac{ G_{- } ( 1 , s) b_0^+ ( s)}
{ G_{+ } ( 1 , s) b_0^- ( s)} \,. \]
For $ l \in \NN $, $ G_+ ( x, l/2) = G_- ( x , l/2 ) $ and 
$ b_0 ^+ ( l/2) \neq 0 $, $ b_0 ^- ( l/2 ) = 0 $. Hence 
$ S ( s ) $ has a simple pole at $ s = l/2 $, $ l \in \NN$.

Another independent solution for $ s = l/2$
is obtained by taking
\[
\begin{split}
&  \partial_s ( G_- ( x, s ) - G_+ ( x , s ) ) |_{ s = l/2} \\
& \ \ = 
x^{ -\frac{l}2 } \left( \sum_{ j = 0 }^{l-1} \partial_s b_j ^- (l/2) x^j 
- 2 b_l ^{-} ( l/2 ) x^l \log x + \Oo ( x^l ) \right) \,, \end{split}
\]
and dividing by 
$ \partial_s b_0^- ( l/2) $ 
gives us a solution with the prescribed leading part, as in 
\eqref{eq:2.scat}:
\begin{equation}
\label{eq:2.4} u_{ l/2} ( x ) = x^{-l/2} + \cdots -2 p_l x^{l/2} \log x
+  \Oo ( x^{l/2} )  \,,
\end{equation}
where $p_l=b_l^-(l/2)/\pa_s b_0^-(l/2)$.
Comparison of the definition of the scattering matrix with the expansions
gives
\begin{equation}
\label{eq:2.5}
p_l = - {\rm{Res}\;}_{ s = l/2  } \; S (s)  \,.
\end{equation}

Comparison of this elementary discussion with the way in which the
invariant operators $ P_k $ are formulated in \S \ref{conf}
provides motivation for Theorem \ref{t:1}.

\section{Scattering on conformally compact manifolds}
\label{scat}

We recall from \S 1 that
the class of asymptotically hyperbolic manifolds, $ ( X , g )$, 
is given by conformally compact manifolds 
\eqref{eq:confc} for which all sectional curvatures tend to $ -1 $ 
as the boundary is approached. That is equivalent to demanding that
\begin{equation}
\label{eq:3.0}
  |d x|_{ \bar g }  = 1  \ \ \text{at $ \ \partial X $} \,, \end{equation}
which depends only on $ g $ -- see \cite{Ma}.
This class of manifolds 
comes with a well developed scattering theory, which originated in  
the study of infinite volume hyperbolic quotients 
$ \Gamma \backslash \HH^{n+1}$ by Patterson, Lax-Phillips, Agmon, 
Guillop\'e, Perry 
and others -- see \cite{Pe}, and references given there. 

The basic facts about the spectrum $\sigma ( \Delta_g )$ of the Laplacian 
of an asymptotically hyperbolic metric were established by 
Mazzeo and 
Mazzeo-Melrose \cite{Ma}, \cite{MM}, \cite{Ma3}.  It is given by 
\[ \sigma ( \Delta_g ) = [ ( n/2)^2 , \infty ) \cup \sigma_{\rm{pp}} 
( \Delta_g ) \,,  \ \  \sigma_{\rm{pp}} ( \Delta_g ) \subset (0, ( n/2)^2 )
\,, \]
where the  pure point spectrum, $ \sigma_{\rm{pp}} ( \Delta_g ) $ (the set
of $L^2$ eigenvalues), is finite. 

More refined statements follow from the main result of \cite{MM} which
is the existence of the meromorphic continuation of  the resolvent 
\[ R ( s ) = ( \Delta_g - s ( n- s ) )^{-1} \,. \]
This traditional choice of the spectral parameter $ s $ corresponds to 
\[ \begin{split} 
\Re s \neq \frac{n}{2} \ & \Longleftrightarrow \ 
 s ( n - s ) \notin [ ( n /2)^2 , \infty )
 \\
\Re s > n \ & \Longrightarrow s ( n - s ) \notin [0, \infty ) \,. 
\end{split}
\]
Since $ \Delta_g \geq 0 $, the second implication and the spectral theorem
show that for $\Re s>n$,
\[ R ( s) : L^2 ( X ) \ \longrightarrow L^2 ( X ) \]
is a holomorphic family of operators. In fact $ R ( s) $ is meromorphic
on $ L^2 ( X) $ for $ \Re s > n/2 $, and on a smaller space continues
meromorphically to $ \CC$:
\begin{prop}(Mazzeo-Melrose \cite[Theorem 7.1, Lemma 6.13]{MM})
\label{p:mm}
The resolvent
\[ R ( s ) \; : \; \DI ( X ) \ \longrightarrow x^s \CI ( X ) 
\]
is meromorphic in $ \{\Re s > n/2\} $, with 
poles having finite rank residues exactly for $ s $ such that $ s ( n - s )  
\in \sigma_{\rm{pp}} ( \Delta_g ) $. It extends meromorphically to
$ \CC $, with poles having finite rank residues, and is regular 
for $\Re s = n/2, \, s\neq n/2$.
\end{prop}
\noindent
In \cite{MM}, the structure of the kernel of $R(s)$ is given in
Theorem 7.1.  The mapping property stated above is proved in Lemma 6.13
only for hyperbolic space, but the proof is valid more generally using only
the properties of $R(s)$ established in Theorem 7.1.  The regularity for
$\Re s = n/2$, $s\neq n/2$ follows from the 
spectral theorem and the absence of embedded eigenvalues -- see 
\cite{Ma3}. 

The line $\Re s = n/2$ corresponds to physical scattering.  
As in the one dimensional case, for such $s$ other than $n/2$, we will
consider generalized eigenfunctions $ u $ of the form
\begin{equation}
\label{expansion}
\begin{split}
&  ( \Delta_g - s ( n - s )) u = 0  \\ 
& u = F x^{n-s} + G x^s \,, \ \ F, G \in \CI ( X ) \,. 
\end{split}
\end{equation}
Note that $n=0$ in the example in \S \ref{ex}, and as there we think
of $ F |_{\partial  X} $ and $ G|_{ \partial X } $ as the
incoming and outgoing scattering data respectively.  
The scattering matrix is the operator which maps $F|_{\pa X}$ to
$G|_{\pa X}$.  For $\Re s > n/2$, it
is more appropriate to think of (\ref{expansion}) as a boundary value
problem.  Then  $F|_{\pa X}$ represents ``Dirichlet data'' and 
$G|_{\pa X}$ ``Neumann data'', so that the analytically continued
scattering matrix can be regarded as a generalized Dirichlet to Neumann
map.  As we shall see, the expansion for $u$ must
in general be modified by the inclusion of a log term if $2s-n\in \N$.    

There are two consequences of Green's identity which are relevant for us.
In the case $\Re s=n/2$, the following identity is known as the {\em
pairing formula}:
\begin{prop}
Let $  \DI ( X ) $ denote functions vanishing to 
infinite order at $ \partial  X $. For $ \Re s = n/2 $ and 
$u_1,u_2$ satisfying
$$
( \Delta_g - s ( n - s ) ) u_i = r_i \in  \DI ( X),
$$
$$
u_i = x^{n-s } f_i + x^{s} g_i + \Oo ( x^{n/2 + 1} ) \,, \ \ 
f_i,g_i \in \CI ( \partial X ), 
$$
we have 
\begin{equation}\label{eq:2.pair}
\int_X \left( u_1 \bar r_2 - r_1 \bar u_2 \right) dv_g  =
( 2 s -n ) \int_{ \partial X } \left( f_1 \bar f_2 - g_1 \bar g_2 
\right) dv_h,
\end{equation}
where, since we chose the defining function $ x $ of $ \partial X $, we
 have a natural choice in the conformal class, $ h=x^2g|_{T\pa X}$. 
\end{prop}

\begin{proof}
This is a standard application 
of Green's identity -- see \cite{Mel} and the proof of Proposition 
\ref{p:witten} below.
\end{proof}

The following analogous identity for $s\in \RR$ will be crucial
in the proof of Theorem~\ref{t:3}.

\begin{prop}
\label{p:witten}
Suppose $s\in \RR$, $s >n/2$, $2s-n\notin \NN$.
Let $u_1$, $u_2$ be real-valued and satisfy (\ref{expansion}) with 
$F_i,G_i \in \CI(X)$.
Then   
\begin{equation}
\label{eq:wit}
 {\rm{pf} } \;
 \int_{ x >\epsilon}[\langle du_1,du_2\rangle -s(n-s)u_1 u_2]
dv_g 
=-n\int_{\pa X} F_1 G_2 dv_h = -n \int_{\pa X} G_1 F_2 dv_h\,,
\end{equation}
where $ {\rm{pf}} $ denotes the finite part of the divergent integral.
\end{prop}
\medskip
\noindent
A special case of \eqref{eq:wit} was discussed by Witten \cite{W}
in a physical context,  for $ X = \HH^{n+1} $, using 
explicit formul{\ae} for $u$ and $G|_{\pa X}$ in terms of $F|_{\pa X}$
given by the Poisson operator and the scattering matrix.  Our
project originated from an attempt to understand that discussion in the
general setting.

\begin{proof}
With $g$ in the form (\ref{normalform}), Green's identity gives 
$$\int_{x>\ep}[\langle du_1,du_2\rangle -s(n-s)u_1 u_2]dv_g
=-\ep^{1-n}\int_{x=\ep} u_1 \pa_x u_2 dv_{h_{\ep}}.$$
Substituting (\ref{expansion}) and
calculating the finite part gives  
$$
\mbox{pf}\int_{x>\ep}[\langle du_1,du_2\rangle -s(n-s)u_1 u_2]
dv_g
=-\int_{\pa X} [s F_1 G_2  + (n-s) F_2 G_1] dv_h.
$$
By symmetry of the left hand side, we deduce that 
$\int_{\pa X} F_1 G_2 dv_h  =\int_{\pa X} F_2 G_1 dv_h$, so the result 
follows. 
\end{proof}

\medskip
\noindent
{\bf Remark.} 
If we assume that $ u_i $'s have expansions \eqref{expansion}
and that
\[  ( \Delta_g - s ( n - s ) ) u_2 = 0 \,, \ \ 
 ( \Delta_g - s ( n - s ) ) u_1 \in \DI ( X) \,, \ \ F_2 = 0 \,, \]
then the same argument gives
\begin{equation}
\label{eq:rem} \int_X u_2( \Delta_g - s ( n -  s ) ) u_1 dv_g = 
(n-2s ) \int_{ \partial X } G_2 F_1 d v_h  \,, \ \ 
s \in \RR \,, \ \ s > n / 2 \,, \end{equation}
which will be useful later.
\medskip

To have the incoming and outgoing data 
$ f_i,g_i $ (or $ F_i|_{\partial X }$, 
$ G_i|_{\partial X } $ in \eqref{expansion})
invariantly defined we introduce density bundles of 
\cite{MelZw}, $ | N^* \partial X|^{s} $, which keep track of changes of the
defining function to first order, or equivalently of the choice of the 
metric in the conformal class\footnote{These bundles are equivalent to the
density bundles of conformal geometry -- see \S \ref{conf}.}:
somewhat informally,
\[ f \in \CI ( \partial X ,  | N^* \partial X|^{s} ) \ \Longleftrightarrow \
f = a |dx|^s \,, \ a \in \CI (  \partial X)  \,. \]
We will often identify
$\CI(\pa X,| N^* \partial X|^{s})$ with $\CI ( \partial X)$ in the presence
of a chosen defining function.

Using Proposition \ref{p:mm} we will construct a
family of
{\em Poisson operators} $\cP(s)$ for $\Re s \geq n/2$, $s\neq n/2$, which
is meromorphic for $\Re s >n/2$ with poles only for
$s$ such that $s(n-s)\in \sigma_{\rm{pp}} ( \Delta_g ) $, and continuous 
up to $\Re s =n/2$, $s\neq n/2$, with the properties: 
\begin{gather}
\label{eq:2.pois}
\begin{gathered}
  \cP ( s ) \; : \; \CI ( \partial X ,  | N^* \partial X|^{n-s} ) 
\; \longrightarrow \CI ( \Xint) \,, 
\\
 ( \Delta_g - s ( n - s ) ) \cP ( s) \equiv 0 \,, \\
 \cP ( s ) f = x^{ n - s } f + o(x^{n-s}) \,, \ \ \mbox{ if } 
\Re s> \frac{n}2,   \\
 \cP ( s ) f = x^{ n - s } f + x^s g + \Oo ( x ^{ n/2 + 1  } ) \,, 
\ \ \mbox {if }\Re s = \frac{n}2, \ \ s \neq \frac{n}{2},
\end{gathered}
\end{gather}
where $g\in \CI(\pa X)$ and we chose trivializations of
the density bundles as indicated above.  If in addition $s\notin
n/2+\NN_0/2$, then we will have
\begin{equation}\label{exptoo}
\cP(s)f \in x^{n-s}\CI(X) + x^s \CI(X)\,.
\end{equation}

As in the example in \S 2, special care will be needed at 
$s\in n/2 + \N/2$.  
To construct $ \cP ( s ) $ for $ \Re s \geq n /2 $ and $ s \notin 
n/2 + \NN_0/2$, we first construct a non-linear 
{\em formal solution} operator:
\begin{equation}
\label{eq:2.formal}
\begin{split}
&  \Phi ( s ) \; : \; \CI ( \partial X)
\; \longrightarrow  x^{n-s}\CI ( X) \,, \\ 
& \text{$ \Gamma ( n+1 - 2 s )^{-1} \Phi ( s) $ 
is entire in $s$,}
\\
& ( \Delta_g - s ( n - s ) ) \Phi ( s)f \in \DI ( X )  \,, \\
& \Phi ( s ) f  = x^{n-s}f +  \Oo (x^{n-s+1}) \,. 
\end{split}
\end{equation} 
This is done by the usual asymptotic construction which we briefly 
recall.  Choose a metric $h$ in the conformal class on the boundary; this
determines a defining function $x$ and product identification near $\pa X$
so that $g$ is given by (\ref{normalform}).  It is easy to see that the
Laplacian, $\Delta_g$, can be written as
\[ \Delta_g = -( x \pa_x )^2 +  n x \pa_x + x^2 \Delta_{h } + x E \,, \;\;
E \in \Diff^2_0 ( X ) \,, \]
that is $E$ can be locally written as a polynomial of degree 2 in $ x
\pa_x $ and $ x \pa_{y_i} $ with coefficients in $\CI(X)$, where $ y_i$ are
local coordinates on $ \partial X $. 

{From} this we see that for $ f_j \in \CI ( \partial X ) $,
\begin{equation}
\label{eq:new}
 ( \Delta_g - s( n -s) ) ( x^{ n - s + j } f_j ) 
= j ( 2 s - n - j ) x^{ n - s + j } f_j +  ( x \Delta_h + {\tilde E} ) 
( x^{ n - s + j + 1 } f_j ) \,, \end{equation}
with ${\tilde E}\in \Diff^2_0 ( X )$.
By an iterative procedure we obtain a formal power series solution 
(see the proof of Theorem~\ref{conformalexpansions} for further details), 
which then gives a solution modulo $ \DI ( X) $ by applying 
Borel's lemma (see for instance \cite[Theorem 1.2.6]{Hor}). It is
here that our operator becomes non-linear: the Borel construction
depends on the sequence of functions $ f_j $.\footnote{The use of Borel's
Lemma could be avoided by truncating the expansions at a large finite
order.  Either direct estimation or a
Banach-Baire type of argument implies that in Proposition \ref{p:mm},
$ \DI ( X ) $ can be replaced by $ x^M \CI ( X) $ for some $ M = M ( s) $,
which is locally constant.  This approach has the advantage that linearity
and continuity are transparent.  However, the use of Borel's Lemma suffices
for our purposes since we are primarily interested in the dependence on $s$
with $f$ fixed.} 
This gives  $\Phi ( s ) $ satisfying \eqref{eq:2.formal}:
\begin{equation}
\label{eq:fexp}
\begin{split}
\Phi ( s) f  =  & fx^{n-s}  +
p_{ 1 , s } f x^{n-s+1} + \cdots \\
& + \; p_{ j , s } f x^{n-s+j}
+  \Oo ( x^{n-s+j + 1} ) \,, 
\end{split} 
\end{equation}
where 
$p _{ j , s }$ is a differential operator on $ \partial X$
of order at most $2[j/2]$, such that 
$$\Gamma(n+1+j-2s)\Gamma(n+1-2s)^{-1}p_{j,s}$$ 
is polynomial in $s$ for each $j$. 

We can now easily prove
\begin{prop}
\label{p:2.poiss}
For $ \Re s \geq n/2 $, $ s \notin n/2 + \NN_0/2 $, $ s ( n - s) \notin
\sigma_{\rm{pp}} ( \Delta_g ) $, there exists a unique linear 
operator $ \cP ( s) $
satisfying \eqref{eq:2.pois}, \eqref{exptoo}.
\end{prop}
\begin{proof}
The existence (of a possibly nonlinear operator) 
comes from the properties of $ \Phi ( s ) $ and Proposition \ref{p:mm}:
\[ 
 \cP ( s ) = ( I - R ( s ) ( \Delta _g - s ( n - s ) ) ) 
 \Phi ( s)  \,. 
\]
To see the uniqueness we first observe that if $ v_j $, $ j = 1,2 $, 
are two solutions in $ x^{ n -s } \CI ( X) + x^ s \CI ( X) $ with the same 
$ x^{n-s} $ leading term, then the formal expansion analysis above shows that
\[ v _1 - v_2 \in x^{ s } \CI ( X)  \,. \]
For $ \Re s > n / 2 $ that would imply that $ v_ 1 - v _2 \in L^2 ( X) $,
contradicting the assumption that $ s ( n - s ) \notin \sigma _{ \rm { pp}} 
( \Delta_g ) $. For $ \Re s = n /2 $ we can apply \eqref{eq:2.pair}
with $ u_1 = v_1 - v_2$ and $ u_2 = R ( n - s ) f $, $ f \in \DI ( X) $, 
which shows that $ v_1 = v_2$.  That $\cP(s)$ is linear and
well-defined on densities follows from uniqueness. 
\end{proof}

Consider the decomposition given above:
\begin{equation}
\label{eq:poiss'}
\cP(s)f = \Phi(s)f - R(s)( \Delta _g - s ( n - s ) ) )  \Phi ( s) f \,, 
\end{equation}
for $s$ as in Proposition~\ref{p:2.poiss}.
We have $\Phi(s)f \in x^{n-s}\CI(X)$ and 
$R(s)( \Delta _g - s ( n - s ) ) )  \Phi ( s) f \in x^s \CI(X)$.  It is
evident from the construction of $\Phi(s)$ that $\Phi(s)f$ may have poles for
$s\in n/2 + \NN/2$.  We shall see that
the contributions from these poles in the two terms cancel one another, so
that $\cP(s)f$ will extend holomorphically across $s\in n/2 +\NN/2$ so long
as $s(n-s)\notin \sigma_{\rm{pp}} ( \Delta_g )$.  

To analyze $ \cP ( s ) $ for  $ s$ near $ n / 2 + \NN / 2 $ we need to 
modify the formal solution operator $  \Phi( s ) $.  Observe that 
$\Phi(n-s)$ satisfies: 
\[
\begin{split}
&  \Phi(n-s) \; : \; \CI ( \partial X)
\; \longrightarrow x^{s } \CI ( X) \,, \\ 
& \text{$\Phi(n-s) $ is  holomorphic in $ \Re s > n /2 $,} 
\\
& ( \Delta_g - s ( n - s ) ) \Phi(n-s)g \in \DI ( X )  \,, \\
& \Phi(n- s ) g = x^{  s } g + {\mathcal O } ( x ^{ s + 1  } ) \,. 
\end{split}
\] 
Fix a metric in the conformal class on $\pa X$ and therefore a defining
function as in \eqref{normalform}.
For $l\in \NN$ and $ s $ near $ n/2 + l/2 $,
consider the modified formal solution operator 
$$\Phi_l(s) \stackrel{\rm{def}}{=} \Phi(s) - \Phi(n-s)p_{l,s}.$$  We certainly have
that $\Phi_l(s)$ is holomorphic for $s$ near $n/2+l/2$, $s\neq n/2+l/2$, and
\[
\begin{split}
&  \Phi_l ( s ) \; : \; \CI ( \partial X)
\; \longrightarrow  x^{n-s}\CI ( X)+x^s \CI(X) \,, \\ 
& ( \Delta_g - s ( n - s ) ) \Phi_l ( s)f \in \DI ( X )  \,, \\
& \Phi_l ( s ) f  = x^{n-s}f +o(x^{n-s}) \,. 
\end{split}
\]
We claim that $\Phi_l(s)$ extends holomorphically across $s=n/2+l/2$
as a map into $\CI(\Xint)$.
Each of $\Phi(s)$ and $p_{l,s}$ has at worst a simple pole, so 
$$ \lim_{ s \rightarrow n/2 + l/2} ( 2s - n - l ) \Phi_l(s)f \,$$ 
certainly exists.  This limit corresponds to 
a formal solution of $ ( \Delta_g - s ( n - s ) ) u =0$  for 
$s=n/2+l/2$, which is of the form $x^{n/2+l/2+1}\CI(X) $, so by 
the iterative formula \eqref{eq:new}, it vanishes identically.

Using 
$$
\lim_{t\rightarrow 0}\frac {x^t - 1}{t} = \log x,
$$
evaluation of the limit as $s\rightarrow n/2+l/2$
shows that
$$\Phi_l(n/2+l/2)f = x^{n/2-l/2}F + x^{n/2+l/2}\log x \,\,G$$
with $F,G\in \CI(X)$, and $F|_{\pa X}=f$, $G|_{\pa X}=-2
\Res_{s=n/2+l/2}p_{l,s}f$. 

For $s$ near $n/2+l/2$ we define 
\begin{equation}
\label{eq:poissl} \cP_{ l } ( s )  \stackrel{\rm{def}}{=} 
( I - R ( s) ( \Delta_g - s ( n - s) ) ) 
 \Phi_l ( s)  \,. \end{equation}
It is clear that if 
$ (n/2)^2 - (l/2)^2 \notin \sigma_{\rm{pp}}(\Delta_g) $, then 
 $\cP_l(s)$ is holomorphic across $s=n/2+l/2$.  
{From} the properties of $\Phi_l(s)$ and Proposition~\ref{p:mm}, it follows
that we have 
\begin{gather*}
 \cP_l ( { n/2 + l/2 } ) \; : \; \CI ( \partial X )   
\ \longrightarrow \ x^{n/2 - l/2}
 \CI ( X) +x^{ n/2 + l/2 } \log x \; \CI ( X) \\ 
 ( \Delta_g - s ( n - s) ) \cP_l ( { n/2 + l/2 } ) = 0 
\,.
 \end{gather*}
For $s\neq n/2+l/2$ we have 
\[ \cP_{ l} ( s) \; : \; \CI ( \partial X ) 
\ \longrightarrow \ \left( x^s \CI ( X) +x^{ n - s} \CI ( X) \right) \cap 
{\rm{ker}} \; ( \Delta_g - s ( n - s) ) \,, \]
so from the uniqueness part of Proposition \ref{p:2.poiss} we see
that 
\[ \cP_{l} ( s) = \cP ( s) \,, \ \ \text{ for $ s $ near $ n/2 + l/2 $, 
\ \ $ s \neq n /2 + l /2 $.} \]
We summarize the discussion of Poisson operators in the following
\begin{prop}
\label{p:2.3}
There is a unique family of Poisson operators 
$$ \cP ( s): \CI ( \partial X ,  | N^* \partial X|^{n-s } )
\longrightarrow \CI(\Xint)$$ 
for $ \Re s \geq n/2, \, s\neq n/2$, which is 
meromorphic in $\{\Re s > n/2 \}$ with poles only for $s$ such that 
$s(n-s)\in \sigma_{\rm{pp} } ( \Delta_g ) $, and continuous up to 
$\{\Re s =n/2\} \setminus \{n/2\}$, such that 
$$(\Delta_g -s(n-s))\cP(s) = 0,$$
with expansions
\[
\begin{array}{ll}
\cP(s)f = x^{n-s}F + x^s G & \mbox{if } s\notin n/2+\NN_0 /2 \\
\cP(s)f = x^{n/2-l/2}F + G x^{n/2+l/2}\log x & \mbox{if } s=n/2+l/2, \,
l\in\NN\,,
\end{array}
\]
for $F,G\in \CI(X)$ such that $F|_{\pa X}=f$.

If $s=n/2 + l/2$, then 
$G|_{\pa X}= -2 p_l f$, where
\[p_l =  \Res _{ s = n /2 + l/2 } p_{ l , s } \,, \]
with the differential operator $ p _{ l , s} $ defined in 
\eqref{eq:fexp}. For $ l = 2k $, 
\[ \sigma_{ 2k } ( p_{2k} ) = c_k \sigma_{2k} ( \Delta_h^{k} ) \,, \] 
with $c_k$ as in \eqref{eq:t.1}.
\end{prop}

The principal symbol of $p_{2k}$ is easily calculated
from the inductive construction of the operators $p_{j,s}$; this will be 
described in more detail in the proof of
Proposition~\ref{conformalexpansions}.  

\medskip
\noindent
{\bf Remark.} A more precise notion of a meromorphic family of 
asymptotic expansions can be given using a Mellin transform in 
$ x $. This gives a holomorphic family of meromorphic
functions with poles corresponding to the exponents in the 
expansions and logarithmic terms to double poles. We do not need 
this precise description here as we will work with the explicit formula
for $ \cP( s ) $.
\medskip

Next we define the scattering matrix 
$$
S(s): \CI(\pa X,| N^* \partial X|^{n-s})\rightarrow
\CI(\pa X,| N^* \partial X|^s)
$$
for $\Re s\geq n/2$, $2s-n\notin \NN_0$, 
$s(n-s)\notin \sigma_{\rm{pp} } ( \Delta_g ) $.  According to 
Proposition~\ref{p:2.3}, for such $s$ 
and for $f\in \CI(\pa X,| N^* \partial X|^{n-s})$
we have $\cP(s)f = x^{n-s}F + x^s G$, so 
we define 
$$S(s)f=G|_{\pa X} \, .$$  The construction of $\cP(s)$ 
in Proposition~\ref{p:2.poiss} gives an explicit realization
of the scattering matrix:  
\begin{equation}\label{scatmatrix}
S(s)f= -\left (x^{-s}R(s)( \Delta _g - s ( n - s ) ) \Phi(s)f\right ) 
|_{\pa X} \,, \ \ 2 s -n \notin \NN_0\,.  
\end{equation}
Using \eqref{eq:2.formal} and Proposition \ref{p:mm}, it follows that 
$S(s)$ is defined for $s$ as indicated above and is meromorphic in $\{\Re s 
>n/2 \}$.  

Applying complex conjugation to \eqref{expansion}
shows that for $ \Re s = n/2 $, $s\neq n/2$, we have
the {\em functional equation}
\[ S( s ) = S ( n - s)^{-1} \,. \]
An application of the pairing
formula \eqref{eq:2.pair} with $ u _j = \cP( s ) f_j $ shows the unitarity 
relation:
\[
\int_{\partial X} \left(   S ( s ) f_1 \overline { S ( s ) f_2 } - 
f_1 \bar f_2 = 0 \right) dv_h \,, \ \ \Re s = \frac{n}2 \,,\ s\neq n/2\,.  
\]

By using the Schwartz reflection principle and then unitarity to establish 
regularity at $s=n/2$, we find a meromorphic extension
of $ S( s) $ to the entire complex plane, regular for $s\in i\RR$, and we
have the following symmetries: 
\[
S( n - \bar s )^* = S ( s )^{-1}  \,, \ \ S ( n -s ) =  S ( s ) ^{-1}  \,.
\] 
For $s\in \RR$, the functional equation and unitarity show that 
\[ S ( s) = S( s )^* \,, \]
which also follows directly from 
Proposition~\ref{p:witten} since $ S ( s ) $ is a real operator for 
$ s $ real.

The construction of the Poisson 
kernels and Proposition \ref{p:2.3} give the following:
\begin{prop}
\label{p:2.5}
The scattering matrix is meromorphic in $ \Re s >  n/2 $ and 
at a pole $ s_0 $ we have 
\[ {\rm{Res}\;}_{ s = s_0   } \; S (s)  = 
\left\{ \begin{array}{ll} { \Pi_{s_0} } & s_0 \notin n/2 + \NN/2 \\
{ \Pi_{n/2 + l/2} } - p_l & s_0 = n/2 + l /2 \,, 
\end{array} \right. \]
where $ p_l$ is as in Proposition \ref{p:2.3} and
\[  \Pi_{s_0} = ( 2 s_0 - n ) \left( x^{-s_0} \frac{1}{2 \pi i } 
\int_{\gamma_{s_0}} 
R ( s ) d s \; x^{-s_0} \right) \big|_{\partial X \times \partial X } \, \]
is non-trivial only if $ s_0 ( n - s_0 ) \in \sigma_{\rm{pp}} ( \Delta_g )
$.  Here $\gamma_{s_0}$ is a small circle in $\CC$ about $s_0$, traversed 
counterclockwise. 
\end{prop}
\begin{proof}
Suppose first that the resolvent is holomorphic near $ n/2 + l/2$. Then
\eqref{eq:poiss'}, \eqref{scatmatrix}, and the holomorphy of $ 
\cP ( s ) $ near $ n/2 + l/2 $ show that
\[ \begin{split}   {\rm{Res}\;}_{ s = n/2 + l/2    }  S ( s )  
& = - \left(  {\rm{Res}\;}_{ s = n/2 + l/2    }  x^{-s}  \Phi ( s) 
\right)  |_{\partial X} \\
& = -   {\rm{Res}\;}_{ s = n/2 + l/2    }   p_{ l , s } = - p _l \,, 
\end{split} \]
where $ p_{ l , s} $ are as in \eqref{eq:fexp}.

When $ s_0 $ is a pole of the resolvent then 
\begin{equation}\label{Rres}
R ( s ) = \frac{ \sum_{ k=1}^{K} \phi_k \otimes \phi_k }{ 
s_0 ( n - s_0 ) - s ( n - s ) } + \widetilde T ( s ) = 
\frac{1}{ 2 s_0 -n } \frac{ \sum_{ k=1}^{K} \phi_k \otimes \phi_k }{ 
s - s_0 } +  T ( s ) \,, 
\end{equation}
where $ T $ and $ \widetilde T $ are holomorphic near $ s_0$, and 
$ \phi_k \in x^{s_0} \CI ( X) \subset 
L^2 ( X) $ are the normalized eigenfunctions, $ ( \Delta_g - 
s_0 ( n - s_0) ) \phi_k = 0 $.  
Let us first consider the case $ s_0 \notin n/2 + \NN/2 $.
We use \eqref{eq:rem} to evaluate
\[ \int _X \phi_k ( \Delta_g - s_0 ( n - s_0) ) \Phi ( s_0) f  dv_g = 
(n -2s_0 ) \int_{\partial X } ( x^{-s_0} \phi_k )|_{\partial X} 
f d v_h \,. \]
Using this and \eqref{Rres} in \eqref{scatmatrix} shows 
that the residue of $ S ( s ) $ at
$ s_0 $ is the operator with kernel
$\sum_{ k=1}^{K} (x^{-s_0}\phi_k)|_{\pa X} \otimes (x^{-s_0}\phi_k)|_{\pa
X}$, which by \eqref{Rres} is $ \Pi _{s_0} $.

If $ s_0 $ is a pole of the resolvent and $ s_0 = n/2 + l/2 $ then 
we replace \eqref{eq:poiss'} by \eqref{eq:poissl} and find the 
residue of $ \cP ( s) $ by the same method as above:
\[  {\rm{Res}\;}_{ s = n/2 + l/2    }  \cP ( s) f = \sum_{ k = 1 } 
^K \phi_k  \int_{\partial X } ( x^{-s_0} \phi_k )|_{\partial X} 
f d v_h \,. \]
Now using this instead of holomorphy of $\cP(s)$ in the argument of the
first paragraph of this proof gives the desired formula for
$   {\rm{Res}\;}_{ s = n/2 + l/2    } S ( s) $. 
\end{proof} 

This shows that the scattering matrix {\em always} has poles of infinite rank 
at $ n/2 + k $, $ k \in \NN $, and possibly at $ n/2 + k - 1/2 $, $ k 
\in \NN $, with the residues given by the operators appearing as the
log term coefficients in the expansions of the Poisson operators.

Uniqueness of $ \cP ( s) $ shows that  $\cP(n)1=1$.  
It then follows from Proposition~\ref{p:2.3} that
$p_n 1=0$ so that $p_{n,s}1$ extends holomorphically across $s=n$.  
Since we always have
$0\notin \sigma_{\rm{pp}} ( \Delta_g ) $, Proposition~\ref{p:2.5} 
shows that $S(s)1$ extends holomorphically across $s=n$.  
Set $S(n)1=\lim_{s\rightarrow n} S(s)1$.  Then we have:

\begin{prop}
\label{p:2.6}
If $ p_{ n, s } $ are defined in \eqref{eq:fexp}, then
\[ S ( n ) 1 = - \lim_{ s \rightarrow n } p_{n,s }1  \,. \]
\end{prop}
\begin{proof}
We have $\cP(s)1 = \Phi(s)1 +x^s G(s)$, where 
the expansion of $\Phi(s)1$ is as in \eqref{eq:fexp} and $G(s)\in \CI(X)$ 
satisfies $G(s)|_{\pa X} = S(s)1$.  Letting $s\rightarrow n$ and recalling
$\cP(n)1 = 1$ shows that $\lim_{s\rightarrow n}(\Phi(s)1 +x^s G(s)) =
1$.  It follows that $\lim _{s\rightarrow n}p_{l,s}1= 0$ for $1\leq l <n$,
and $\lim_{s\rightarrow n} (p_{n,s}1 + S(s)1)=0$ as desired.
\end{proof}

\begin{figure}[htb]
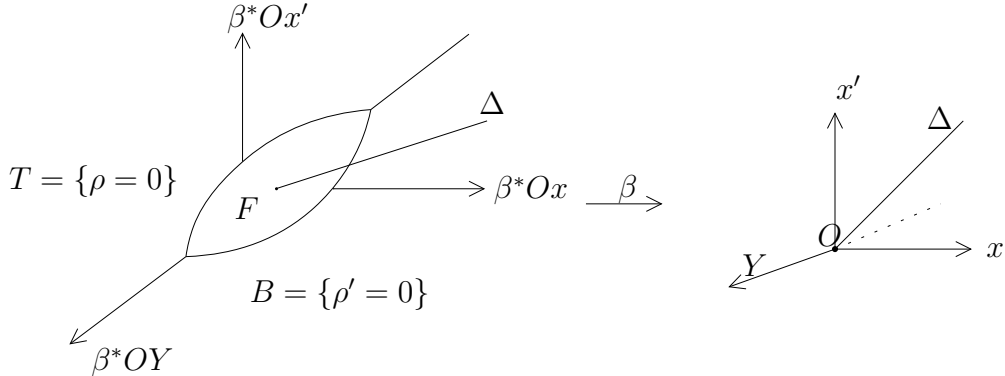

\placerfig{sa1}{15}
\caption{The boundary faces of the blown-up space $ X \times_0 X $.
The variable $ Y $ stands for the defining function of the diagonal
in $ \partial X \times \partial X $, $ Y = y-y' $, in local 
coordinates.}
\label{fig:bup} 
\end{figure}

We next discuss the distributional kernels of the operators $R(s)$, $\cP(s)$,
and $S(s)$.  
The structure of the kernel of $ R ( s ) $ is best understood
on a resolved product space. Since it is not 
central to our discussion, we will only sketch the construction.
We first recall the blow-down map
of \cite[Sect.3]{MM}:
\[ \beta \; : \; X \times_0 X \ \longrightarrow \ X \times X \,, \]
where $X \times_0 X $ is the blow-up of $X\times X$ along the boundary
diagonal, 
illustrated in Fig.\ref{fig:bup}. The restriction to the boundary 
gives us a blow-up of the diagonal in $ \partial X \times \partial X $
with the corresponding blow-down map:
\[ \beta_\partial  \; : \; T \cap B \ \longrightarrow \ \partial X
\times \partial X \,, \]
where $ T $ and $ B $ are the top and bottom faces in the blow-up, 
as shown in Fig.\ref{fig:bup}. 
Let 
$ \rho $, $ \rho' $ be the defining functions of the top and bottom 
faces respectively. Then we have 

\begin{prop}
(Mazzeo-Melrose \cite[Proposition 6.2, Theorem 7.1]{MM})
\label{p:mm2}
The kernel of the resolvent on the resolved space $ X \times_0 X $,
$ \beta^* R ( s ) $, and away from its poles, can be written as 
\[  \beta^* R ( s) = R' ( s ) + R'' ( s) \,, \]
where $ R' ( s ) $ is supported {\em away} from $ T \cup B $, and
\[ R'' ( s ) = ( \rho)^s ( \rho')^s K ( s) \,, \ \ 
K ( s) \in \CI ( X \times_0 X ) \,. \]
\end{prop}

The pairing formula \eqref{eq:2.pair} gives an expression for the kernel of
$ \cP ( s) $
in terms of that of the resolvent, much in the same spirit as in the
classical  derivation of the Poisson kernel from the Green function:
\begin{prop}
\label{p:2.4}
For $ f \in \CI ( \partial X ) $ and $\Re s \geq n/2$, $s\neq n/2$, we have
a formal expression 
\begin{equation}
\label{eq:2.poiss}
 \cP ( s )f  = ( 2s -n ) ( R ( s) {x}^{-s} )  ( f \otimes \delta_0 ( x ) )
 \, ,  
\end{equation}
where the right hand side is understood as the restriction to the boundary
of the distributional kernel of $ R ( s ) x^{-s} $ paired with $ f $. 
\end{prop}
\begin{proof}
First let $\Re s = n/2$.  For any  $ g \in \DI ( X) $  and $ f
\in \CI ( \partial X )  $ apply \eqref{eq:2.pair} with $ u_1 = \cP ( s ) f$ 
and $ u_2 = R ( n- s ) g $. We also note that for $ \Re s = n/2 $ we
have 
\[   R ( n -s )^* = R ( s )  \,. \]
That gives
\[ 
\begin{split}
\int_X \cP ( s ) f \bar g & = ( 2 s -n ) \int_{\partial X} \overline{( 
x^{- n + s} R (n- s) 
g )} |_{\partial X }  f \\
& = ( 2 s -n ) \int_X \left((  R ( s) x^{-s} ) 
( f \otimes \delta_0 ( x) ) \right) \bar g \,, \end{split} \]
proving the claim. That the restriction of $ R ( s) x^{-s} $ is well 
defined as a distribution on the boundary follows from Proposition 
\ref{p:mm2}.  The result extends to $\Re s > n/2$ by analytic
continuation.  
\end{proof}
\noindent
We now continue $\cP(s)$ meromorphically to $\CC$ so that
\eqref{eq:2.poiss} holds for $\Re s< n/2$ as well.

Using the structure of the resolvent described in Proposition \ref{p:mm2}
we can give a precise description of the kernel of $ S ( s ) $ as 
a pseudodifferential operator: away from its poles, we have
\begin{gather}
\label{eq:pseud}
\begin{gathered}
  S (s ) \in \Psi^{ 2 \Re s -n } ( \partial X  ; 
 | N^* \partial X|^{n- s} ,  | N^* \partial X|^{ s}  )  \,, \\
\sigma ( S ( s) ) = 2^{ n - 2s } \frac{ \Gamma( n/2 - s) }{ \Gamma ( 
s - n /2 ) }\sigma ( \Delta_h ^{s - n/2 } ) .
\end{gathered}
\end{gather}

To see this we follow \cite{GZ} and obtain the following formula
for the distributional kernel of the scattering matrix in terms of
the resolvent on the resolved product space:
\begin{equation}
\label{eq:2.scatt}
S( s) = ( 2s -n ) 
(\beta_\partial)_* \left( \beta^* \left( x^{ -s } {x'}^{-s} 
R ( s) \right) \big|_{T \cap B} \right) \,, \end{equation}
(see also \cite[Proposition 4.4]{JS}).
This is derived using \eqref{eq:poiss'}, \eqref{scatmatrix} 
and Proposition \ref{p:2.4}, which give 
\[ \begin{split}
S ( s) f &  = ( x^{ -s } \cP ( s) f   - x^{- s }\Phi(s) f )|_{\partial X} 
\\
& = \left(( 2s -n  )  x^{-s } ( R ( s) x^{-s} ) ( f \otimes \delta_0 ( x) ) 
- x^{- s } \Phi(s)f \right)|_{\partial X } \,. 
\end{split} \]
For $ \Re s < n/2 $ and $ s $ away from poles of $ R ( s) $ we see that
\[ S ( s) f =
 \left(( 2s -n  )  x^{-s } ( R ( s) x^{-s} ) ( f \otimes \delta_0 ( x) ) 
 \right)|_{\partial X } \,.\]
Proposition \ref{p:mm2} shows that
for $ \Re s \ll 0 $, and $ s $ not a pole, we have
\[ \left( x^{-s } ( \beta_{* } R ' ( s) x^{-s} ) ( f \otimes \delta_0 ( x) ) 
 \right)|_{\partial X } \equiv 0 \,, \]
and 
\[ \left( x^{-s } ( \beta_{* } R '' ( s) x^{-s} ) ( f \otimes \delta_0 ( x) ) 
 \right)|_{\partial X }  = 
\left( x^{-s } R ( s) x^{-s} ) ( f \otimes \delta_0 ( x) ) 
 \right)|_{\partial X }  \,, \] 
and this gives \eqref{eq:2.scatt} for $\Re s \ll 0$.  Strictly speaking we
need more information about the $ R' ( s) $ term than stated in Proposition
\ref{p:mm2} but those properties are stated in \cite[Theorem 7.1]{MM}.
We shall see below that \eqref{eq:2.scatt} continues meromorphically in
$s$.  

The fact that the scattering matrix is a pseudodifferential operator of
order $ 2 \Re s - n $ follows from \eqref{eq:2.scatt}. As shown in 
Fig.\ref{fig:bup}, near $ T \cap B $ we can use $y$, $ |Y|$, $ \omega = 
Y/ |Y|$, and $ \rho, \rho' $, the defining functions of $ T $ and $ B $, as
coordinates, where $Y=y-y'$.  In terms of these functions, 
\[ x = |Y| \rho \,, \ \ x' = |Y| \rho' \,, \]
and hence, using Proposition \ref{p:mm2}, 
\begin{equation}
\label{eq:cont}
 S ( s) = (2s -n )  
(\beta_{\partial})_* \left( |Y|^{ -2s } K ( s, y, |Y| , \omega, ) 
\right) \,, \ \ 
K(s,\cdot) \in \CI ( [ \partial X \times \partial X ; 
\Delta_{\partial X } ] )\,, 
\end{equation}
where $  [ \partial X \times \partial X ; \Delta_{\partial X } ] $ is
$ \partial X \times \partial X $ blown up along the diagonal, that
is the space where instead of $ y, y' \in \partial X $ as  coordinates,
we use $ y, | Y | , \omega $. 

Hence for $ \Re s \ll 0 $ the distributional 
kernel of $ S ( s) $ is a distribution conormal to the diagonal and hence 
$ S( s )$ is a pseudo-differential operator. We then continue 
$ |Y|^{-2s} $ meromorphically in $ s $ following \cite[Theorem 3.2.4]{Hor},
and that gives the first part of \eqref{eq:pseud}. The parametrix
construction in \cite{MM} gives $ K $, and leads to a computation of 
the principal symbol -- see \cite{JS}.

\section{Conformal Geometry}
\label{conf}
In this section we show how the conformally invariant powers of
the Laplacian may be derived from the formal Poincar\'e metric associated
to the conformal structure, review the definition and basic properties
of Branson's $Q$-curvature, and remark on the interpretation of the
scattering matrix for 
a Poincar\'e metric as a family of conformally invariant
pseudo-differential operators. 

The invariant powers of the Laplacian act on conformal densities.  The
metric bundle of a conformal manifold $(M,[h])$ is the ray subbundle 
${\mathcal G}\subset S^2T^*M$ of multiples of the metric:  if $h$ is a
representative metric, then the fiber of ${\mathcal G}$ over $p\in M$ is
$\{t^2h(p):t>0\}$. The space of conformal densities of weight $w\in \CC$ is 
\[ \E ( w ) = \CI ( M ; {\mathcal G}^{-\frac{w}2} ) \,, \]
where by abuse of notation we have denoted by ${\mathcal G}$ also the line
bundle associated to the ray bundle defined above.  
A choice of representative $h$ for the conformal structure induces an
identification $\E (w)\simeq \CI (M)$; if $\hat{h}=e^{2\up}h$ then the 
corresponding elements of $\CI (M)$ transform by $\hat{f}=e^{w\up}f$. 
If $\partial X = M$ as in the
previous section, a conformally compact metric on $X$ with conformal
infinity $[h]$ determines an isomorphism 
\[ \E ( w ) \simeq \CI ( M , | N^* \partial X |^{-w} ) \,. \]
The invariance property \eqref{eq:1.inv} can be reformulated as the
statement that $P_k$ is an invariantly defined operator
$$  P_k : \E ( -  n/2 + k )\rightarrow \E ( -  n/2 - k ).$$
The main result of \cite{GJMS} is the following existence theorem.
\begin{prop}
Let $k\in \N$ and $k\leq n/2$ if $n$ is even.  There is a 
conformally invariant natural differential operator
$P_k:\E(-n/2 + k)\rightarrow \E(-n/2-k)$ 
with principal part equal to that of $\Delta^k$.  
\end{prop}
\noindent
Such operators are not unique.  The construction in \cite{GJMS} is in
terms of the ambient metric of \cite{FG}.  We shall first give here a 
different construction of invariant operators based on the Poincar\'e
metric and then show that the two constructions give the same operators. 

We begin by recalling from \cite{FG} the formal Poincar\'e
metric associated to a conformal structure $(M,[h])$.  Given a
representative metric $h$, one considers metrics on 
$M\times [0,1]$ of the form \eqref{normalform}.
The Einstein equation $\Ric(g)+ng=0$ can be calculated directly in terms of 
$h_x$ and the formal asymptotics of solutions studied; see
\cite{G}.  If $n$ is odd, then there is a unique formal smooth solution
$h_x$ to 
$$\Ric(g)+ng= \Oo (x^{\infty})$$ 
which is even in $x$.  If $n$ is even,  
the condition $\Ric(g)+ng= \Oo (x^{n-2})$ uniquely determines 
$h_x \mod  \Oo (x^n)$, which is even in $x$ (mod $ \Oo (x^n)$).  Although in
general smooth solutions do not exist to higher orders, 
the condition 
$$\tr_g(\Ric(g)+ng)= \Oo (x^{n+2})$$
can be satisfied and uniquely 
determines the $h$-trace of the $x^n$ coefficient in $h_x$; this is the
vanishing trace condition referred to in the introduction.  The 
indicated Taylor coefficients of $h_x$ are determined inductively from
the equation and are given by polynomial formulae in terms of $h$, its
inverse, and its curvature tensor and covariant derivatives thereof.

Since any asymptotically hyperbolic metric can be put uniquely in the form
\eqref{normalform} upon choosing $h$, it follows that the equivalence class
of the solution $g$ up to 
diffeomorphism and up to terms vanishing to the indicated orders is
uniquely determined by the conformal structure.  This equivalence class is
called the formal Poincar\'e metric associated to $[h]$.  When $n$ is even,
the higher order terms in $h_x$ are not determined; however for
simplicity in statements below, we shall
restrict consideration to $h_x$ which are smooth and even in $x$ to all
orders.  If $X$ is a manifold with $\partial X=M$, any metric on $X$ whose
restriction to a collar neighborhood of $M$ is in this equivalence class is
called a Poincar\'e metric associated to $[h]$.  

Let $g$ be a Poincar\'e metric and $h$ a
representative for the conformal infinity. If $x$ is a defining 
function such that $x^2 g|_{TM}=h$, and $f\in \CI(M)$ represents a section
of $\E (w)$, then it is a conformally invariant statement to require that a
function $u$ on $X$ be asymptotic to $x^{-w} f$.  The invariant operators
$P_k$ arise from solving
\begin{equation}\label{infiniteorder}
(\Delta_g - s(n-s))u= \Oo (x^{\infty})
\end{equation}
for $u$ with such asymptotic behaviour.  The characteristic
exponents of $\Delta_g - s(n-s)$ are $s, n-s$, so generically 
solutions behave like $x^s$ and $x^{n-s}$.  As in the previous section,
we take $u$ 
asymptotic to $x^{n-s}f$, which according to the above remarks means that
$f$ is to be interpreted as a density of weight $w=s-n$.  

We saw in \S 3 that if $g$ is asymptotically hyperbolic and $s-n/2\notin
\N/2$, then for any $f\in \CI(M)$ there is a formal solution  
$u\mod  \Oo (x^{\infty})$ to (\ref{infiniteorder}) of the form $u=x^{n-s}F$
with $F\in \CI(X)$ and $F|_M=f$.  As we shall see below, for Poincar\'e
metrics this holds 
if only $s-n/2\notin \N$.  In \S 3, solutions for $s-n/2\in \N/2$ were
constructed as a limit of solutions for nearby $s$.  Here we make a direct
analysis at the exceptional values of $s$ and obtain the invariant
operators as obstructions to the existence of formal smooth solutions. 

\begin{prop}\label{conformalexpansions}
Let $(X,g)$ be a Poincar\'e metric associated to $(M,[h])$ and let 
$f\in \CI (M)$.  If $k\in \N$, there is a formal solution of
(\ref{infiniteorder}) for $s=n/2+k$
of the form  
\begin{equation}\label{logform}
u=x^{n/2-k}(F + G x^{2k} \log x)
\end{equation}
with $F,G\in \CI(X)$ and with $F|_M=f$.
$F$ is uniquely determined $\mod \Oo  (x^{2k})$ and $G$ is
uniquely determined $\mod \Oo (x^{\infty})$.  Moreover, 
\begin{equation}\label{logterm}
G|_M = -2c_k P_kf,
\end{equation}
where $P_k$ is a differential operator on $M$ with principal part 
$\Delta^k$.

If $n$ is odd and $k\in \N$ or if $n$ is even and $k\leq n/2$, then $P_k$
depends only on $h$ and defines a conformally invariant operator
$:\E(-n/2 + k)\rightarrow \E(-n/2-k)$.
\end{prop}

\begin{proof}
Choose a representative metric $h$ and write $g$ in the form
(\ref{normalform}).  
A straightforward calculation shows that 
$[\Delta_g-s(n-s)]\circ x^{n-s}= x^{n-s+1}{\mathcal D}_s$, where 
\begin{equation}\label{Laplace}
{\mathcal D}_s=-x\pa_x^2 +(2s-n-1-\frac{x}{2} h^{ij}h_{ij}')\pa_x
-\frac{n-s}{2}h^{ij}h_{ij}'+x\Delta_{h_x}.
\end{equation}
Here $h_{ij}$ denotes the metric $h_x$ with $x$ fixed, and 
$h_{ij}'=\pa_x h_{ij}$.
Therefore for $f_j\in \CI(M)$, one has
\begin{equation}\label{iterate}
{\mathcal D}_s(f_jx^j)=j(2s-n-j)f_jx^{j-1} + \Oo(x^j).
\end{equation}

If $2s-n\notin \N$, then we can use (\ref{iterate}) to construct a smooth
solution of (\ref{infiniteorder}) inductively.  Beginning with $f_0=F_0=f$, 
define $f_j, F_j$ for $j\geq 1$ by
\begin{equation}\label{fj}
\begin{split}
& j(2s-n-j)f_j=-(x^{1-j}{\mathcal D}_s(F_{j-1}))|_{x=0} \\
& F_j=F_{j-1}+f_jx^j.
\end{split}
\end{equation}
By (\ref{iterate}) we have ${\mathcal D}_sF_j=\Oo(x^j)$ so that the 
definition of $f_j$ makes sense.  Then $F=\sum f_jx^j$ is a
formal solution of (\ref{infiniteorder}).  Observe that since $h_x$ is even
in $x$, ${\mathcal D}_s$ maps even functions to odd and vice versa.
Therefore  $f_j=0$ for $j$ odd.  For $j=2k$ even, an easy induction 
shows that $f_{2k}$ takes the form 
\begin{equation}\label{defP}
f_{2k}=c_{k,s}P_{k,s}f, \qquad
c_{k,s}=(-1)^k\frac{\G(s-n/2-k)}{2^{2k}\,k!\,\G(s-n/2)},
\end{equation}
where $P_{k,s}$ is a 
differential operator on $M$, the principal part of which
agrees with that of $\Delta_h^k$.
Since the 
Taylor expansion of $h_x$ is determined (to the appropriate order for
$n$ even) in terms of $h$, one sees by counting derivatives 
that if $k\in \N$ and $k\leq n/2$ if $n$ is even,
then $P_{k,s}$ also depends only on $h$ and is a natural
differential operator with coefficients which are polynomial in $s$.   

If $2s-n=l\in \N$, then the corresponding coefficient $2s-n-j$ vanishes for
$j=l$ in the first equation of (\ref{fj}).  For $l$ odd, by parity
considerations 
it follows that the right hand side of this equation also
vanishes, so $f_l$ can be chosen arbitrarily (for example $f_l=0$ to
preserve parity) and the induction continued to infinite order.

However, if $2s-n=2k$ is even, then the right hand side of the first
equation of (\ref{fj}) need not vanish for $j=2k$ and there is an
obstruction to solving with $F$ smooth.  This is of course reflected in the
pole of  $c_{k,s}$ at $s=n/2+k$.  This obstruction can be incorporated into
a log term and the formal solution continued to higher order as follows.
Observe that  
\begin{equation}\label{log}
\begin{split}
{\mathcal D}_s(g_jx^j\log x)= & 
j(2s-n-j)g_jx^{j-1}\log x + (2s-n-2j)g_jx^{j-1} \\
& \ \ + \; \Oo(x^j\log x) \,. \end{split} 
\end{equation}
Therefore if we take
\begin{equation}\label{g2k}
g_{2k}=(2k)^{-1}(x^{1-2k}{\mathcal D}_s(F_{2k-1}))|_{x=0},
\end{equation}
choose $f_{2k}$ arbitrarily,
and set 
$$F_{2k}=F_{2k-1}+g_{2k}x^{2k}\log x +f_{2k}x^{2k} \,,$$ 
then we have
${\mathcal D}_sF_{2k}= \Oo (x^{2k}\log x)$.  Using (\ref{iterate}) and 
(\ref{log}), it is easily seen that the construction can be continued to
all higher orders to obtain a solution of the form (\ref{logform}) as
claimed.  We have $G|_{x=0}=g_{2k}$, which is given by a differential
operator on $M$ of the claimed form by the same reasoning as above.  In
fact, one sees easily that (\ref{logterm}) holds with $P_k=P_{k,n/2+k}$ and  
$$c_k=\Res_{s=n/2+k}c_{k,s}=(-1)^k[2^{2k}k!(k-1)!]^{-1} \,. $$

Suppose now we change the conformal representative $h$.  We obtain a
different defining function $x$, a different product identification for
$X$, and a different representation (\ref{normalform}).  However, $g$
and $u$ remain unchanged, so by  
uniqueness we deduce that $G|_{x=0}$ must transform as a density, proving
the conformal invariance of $P_k$.
\end{proof}

\begin{rem}
Equation (\ref{Laplace}) holds also for general asymptotically hyperbolic
metrics in the form \eqref{normalform} and can be used to explicitly
compute the operators $p_{j,s}$ of 
\S 3.  Of course in general it no longer need be the case that $h_x$ is
even in $x$.  Therefore log terms may also occur for $2s-n$ odd, and
the scattering matrix may have poles for such $s$ according to
Proposition~\ref{p:2.5}. 
\end{rem}

Note that if $n$ is even, then 
$P_{n/2}$ is invariantly defined from $\E(0)$ to $\E(-n)$, that is,
from $\CI(M)$ to the space of volume densities.
Since the constant function $1\in \E(0)$ has a smooth extension annihilated
by $\Delta_g$, we have $P_{n/2} 1 =0$.  Therefore $P_{n/2}$ has zero
constant term. 

We next recall from \cite{GJMS} the original construction of the invariant
operators via the ambient metric.  Denote by 
$\pi:\cG\rightarrow M$ the natural projection of the metric
bundle, and by ${\bf g}$ the tautological  
symmetric 2-tensor on $\cG$ defined for $(p,h)\in \cG$
and $X,Y\in T_{(p,h)}{\mathcal G}$ 
by ${\bf g}(X,Y)=h(\pi_*X,\pi_*Y)$.  There are dilations 
$\delta_s:{\mathcal G}\rightarrow {\mathcal G}$ for $s>0$ given by
$\delta_s(p,h)=(p,s^2h)$, and we have  
$\delta_s^*{\bf g} = s^2{\bf g}$.  Denote by $T$ the infinitesimal dilation
vector field $T=\frac{d}{ds} \delta_s |_{s=1}$.  Define the ambient space 
${\tilde {\mathcal G}}={\mathcal G}\times (-1,1)$.  Identify ${\mathcal G}$
with its image under 
the inclusion $\iota:{\mathcal G}\rightarrow {\tilde {\mathcal G}}$ given
by $\iota(h)=(h,0)$ for $h\in {\mathcal G}$.  The dilations $\delta_s$ and
infinitesimal generator $T$ extend naturally to ${\tilde {\mathcal G}}$.  

The ambient metric ${\tilde g}$ is a Lorentzian metric on ${\tilde \cG}$
which satisfies the initial condition
$\iota^* {\tilde g}={\bf g}$, is homogeneous in the sense that 
$\delta_s^* {\tilde g}=s^2{\tilde g}$, and is an asymptotic solution of 
$\Ric({\tilde g})=0$ along $\cG$.  For $n$ odd, these conditions uniquely
determine a formal power series expansion for ${\tilde g}$ up to
diffeomorphism, but for $n$ even and $n>2$, a formal power series solution
exists in general only to order $n/2$.   

An element of $\E(w)$ can be regarded as a homogeneous function of degree
$w$ on $\cG$.  One of the ways that
$P_k$ is derived in \cite{GJMS} is as the obstruction to extending 
${\tilde f}\in \E(-n/2+k)$, regarded as such a homogeneous function, to a   
smooth function ${\tilde F}$ on $\tilde \cG$, such that ${\tilde F}$ is
also homogeneous of degree $-n/2+k$ and satisfies 
\begin{equation}\label{harmonic}
\tilde \Delta {\tilde F}= \Oo (\rho^{\infty}),
\end{equation} 
where $\tilde \Delta $ denotes
the Laplacian in the metric ${\tilde h}$.  Similarly to
Proposition~\ref{conformalexpansions}, the Taylor expansion of ${\tilde F}$
is formally determined to order $k-1$ in $\rho$, but there is an
obstruction at order $k$ which defines the operator $P_k$.  

\begin{prop}
Let $k\in \N$ with $k\leq n/2$ if $n$ is even.
The conformally invariant operator $P_k$ defined in \cite{GJMS} in terms
of the ambient metric agrees with the operator of
Proposition~\ref{conformalexpansions}. 
\end{prop}

\begin{proof}
As described in \cite{FG}, the formal
Poincar\'e metric associated to a conformal structure can be 
constructed from the ambient metric and vice versa; the two constructions
are equivalent.  We review this equivalence.

In the ambient space ${\tilde \cG}$, the equation ${\tilde g}(T,T)=-1$
defines a hypersurface 
$\Xint$ which lies on one side of $\cG$ and which intersects
exactly once each dilation orbit on this side of $\cG$.  
The Poincar\'e metric $g$ is the pullback to $\Xint$ of 
${\tilde g}$.  The equation $\Ric({\tilde g})=0$ is equivalent
to $\Ric(g)=-ng$.  To see this, one uses the normal form from \cite{FG} for
ambient metrics:  in
suitable coordinates on ${\tilde \cG}$, the ambient metric takes the form
\begin{equation}\label{ambient}
{\tilde g}= 2tdtd\rho + 2\rho dt^2 +t^2{\overline h}_{\rho}.
\end{equation}
Here $\rho$ is a defining function for $\cG\subset {\tilde \cG}$, $t$ is
homogeneous of degree 1 with respect to the dilations on ${\tilde \cG}$,
and ${\overline h}_{\rho}$ is a smooth 1-parameter family of metrics on
$M$. 
In these coordinates
we have $T=t\pa_t$, so $\Xint=\{2\rho t^2=-1\}$.  Introduce a 
new variable $x=\sqrt{-2\rho}$ and set $s=xt$ so that 
$\Xint=\{s=1\}$.  A straightforward calculation shows that (\ref{ambient})
becomes 
\begin{equation}\label{ambient2}
{\tilde g}=s^2g -ds^2,
\end{equation}
where $g$ is given by (\ref{normalform}) with $h_x={\overline h}_{\rho}$.
The equivalence of $\Ric({\tilde g})=0$ and $\Ric(g)=-ng$ is a
straightforward calculation given the relationship (\ref{ambient2}) 
(see Proposition 5.1 of \cite{GL}).  Note that 
$h_x$ is automatically even in $x$.  

The equations (\ref{infiniteorder}) and (\ref{harmonic}) are equivalent.
To see this, rewrite (\ref{ambient2}) as 
$${\tilde g}=s^2(g -ds^2/s^2)\,, $$ 
and 
transform under a conformal change
to obtain 
$${\tilde \Delta} = s^{-2}[\Delta_g +(s\pa_s)^2+ns\pa_s] \,. $$
If ${\tilde F}$ is homogeneous of degree $w$, we therefore have
$${\tilde \Delta} {\tilde F}= s^{-2}[\Delta_g +w(w+n)]{\tilde F},$$
so \eqref{harmonic} is equivalent to (\ref{infiniteorder}) for $s=n+w$
and  
$u={\tilde F}|_{\Xint}$.  One may recover ${\tilde F}$ from $u$ via
homogeneity by ${\tilde F}=s^w u = t^w x^w u$.  
In order for ${\tilde F}$ to be smooth up to $\rho=0$, we  
require therefore that $x^w u$ be smooth up to $x=0$ (and be even in $x$).
Thus the two extension problems are equivalent, so the normalized
obstruction operators must agree.
\end{proof}

We next define the $Q$-curvature as in \cite{B}.  For this discussion
we shall denote by $P^n_k$ the operator $P_k$ in dimension $n$.  
Fix $k \in \NN$.  One
consequence of the construction of Proposition~\ref{conformalexpansions} 
is that the operator $P^n_k$
is natural in the strong sense that $P^n_k f$ may be written as a linear
combination of complete contractions of products of covariant derivatives 
of the curvature tensor of a representative for the conformal structure
with covariant derivatives of $f$, with coefficients which are rational in
the dimension $n$.  Also, it follows from the fact that the zeroth order
term of ${\mathcal D}_s$ in (\ref{Laplace}) has a factor of $n-s$, that 
the zeroth order term of $P^n_k$ may
be written as $(n/2-k) Q^n_k$ for a scalar Riemannian invariant $Q^n_k$  
with coefficients which are rational in $n$ and regular at $n=2k$.  (This
is of course consistent with the fact mentioned above that $P^n_{n/2} 1
=0$.)  The $Q$-curvature in even dimension $n$ is then defined as 
$Q=Q^n_{n/2}$.  

We may also consider the $Q$-curvature as arising in a similar way from the
zeroth order terms of the operators $P_{k,s}$ with $n$ fixed but as $s$
varies.  For the same reason as above, the zeroth
order term of $P_{k,s}$ is of the form 
\begin{equation}\label{Qs}
P_{k,s}1=(n-s)Q_{k,s}
\end{equation}
for a scalar Riemannian invariant $Q_{k,s}$ which is polynomial in $s$.  
Taking $s=n/2+k$ and recalling that $P_k=P_{k,n/2+k}$ shows that
$Q^n_k=Q_{k,n/2+k}$.  In particular, 
\begin{equation}\label{Qform}
Q=Q_{n/2,n}.
\end{equation}

If $g$ is a Poincar\'e metric with conformal infinity $[h]$, then according
to \eqref{eq:pseud}, the scattering
matrix $S(s)$ is a family of pseudodifferential operators on $M$, which is
conformally invariant in the sense that it acts 
invariantly on conformal densities.  Even though $S(s)$ depends on the
choice of Poincar\'e metric $g$, it is shown in \cite{JS} that its full
symbol depends only on the 
infinite jet of $g$ at $\pa X$, which, as discussed above, is determined by
$[h]$ (to the appropriate order for $n$ even).  
One can therefore view the symbol as determined by the conformal
structure, and the choice of $g$ as a
geometrically natural means of fixing smoothing terms to obtain
a globally well-defined operator having this symbol.
Peterson \cite{Pet} has defined an analogous family of symbols directly by
analytic continuation from the differential operators $P_k$.

Branson \cite{B} derived the transformation law \eqref{eq:1.q} for $Q$ by
analytic continuation in the dimension.  It also follows easily from
Theorems 1 and 2 by analytic continuation in $s$.  The conformal invariance
of $S(s)$ is equivalent to $\widehat{S(s)} = e^{-s\up}S(s)e^{(n-s)\up}$, 
where here $S(s)$ is realized as an operator on $\CI(M)$ corresponding to
the choice of conformal representative.  Therefore,
$$
e^{s\up}\widehat{S(s)}1 = S(s)1 + S(s)(e^{(n-s)\up}-1).
$$
Letting $s\rightarrow n$ and applying Theorems 1 and 2 immediately yields
\eqref{eq:1.q}. 

\section{Proofs of the main results}
\label{pr}

Theorems 1 and 2 follow from the corresponding results of \S 3 
in the special case when $g$
is a Poincar\'e metric associated to the conformal structure.  It follows
from \eqref{defP} that for
Poincar\'e metrics, the differential operators $p_{l,s}$ of \eqref{eq:fexp}
are given by
\begin{equation}\label{corr1}
p_{l,s}=0 \mbox{ for $l$ odd, } \qquad p_{2k,s}=c_{k,s}P_{k,s},\,\,
k\in \NN.
\end{equation}
Thus for the residues we obtain 
\begin{equation}\label{corr2}
p_l=0 \mbox{ for $l$ odd, } \qquad p_{2k}=c_k P_k,\,\, k\in \NN.
\end{equation}

\medskip
\noindent
{\em Proof of Theorem 1.}
This is immediate from Proposition \ref{p:2.5} and \eqref{corr2}. 
\stopthm

\noindent
Observe that Proposition \ref{p:2.5} and \eqref{corr2} also show that for
Poincar\'e metrics, $S(s)$ has no pole at $s=n/2+l/2$ if $l\in \NN$ is
odd. 

\medskip
\noindent
{\em Proof of Theorem 2.}
We apply Proposition \ref{p:2.6}.  By \eqref{corr1} and
\eqref{Qs}, \eqref{Qform}, we obtain 
$$
S(n)1=-\lim_{s\rightarrow n}c_{n/2,s}P_{n/2,s}1=\Res_{s=n}c_{n/2,s}\, Q= 
c_{n/2} Q.
$$
\stopthm

\noindent
We note that Proposition \ref{p:2.6} and \eqref{corr1}
also show that for Poincar\'e metrics, $S(n)1=0$ if $n$ is odd. 

\medskip
\noindent
{\em Proof of Theorem 3.}
Let $n$ be even and let $g$ be a Poincar\'e metric.
We first recall (as described in \cite{G}) the expression for $L$ in terms
of the expansion of the volume form of $g$.
If $h$ is a representative for the conformal infinity of $g$, we may write
$g$ in the form (\ref{normalform}), and 
from the fact that $g$ is a Poincar\'e metric it follows that 
\begin{equation}\label{volumeform}
dv_g =x^{-n-1}(1+v^{(2)}x^2+(\mbox{even powers})+v^{(n)}x^{n}+ \ldots)
dv_h dx,
\end{equation}
where each $v^{(2j)}$ for $1\leq j\leq n/2$ is a smooth function on $M$
determined by $h$.  Integration yields \eqref{eq:1.L}, and shows that 
\begin{equation}\label{Lform}
L=\int_M v^{(n)}dv_h.
\end{equation}

Now take $s$ near but not equal to $n$ and
apply Proposition~\ref{p:witten} with $u_1=u_2=\cP(s)1$, which we
now denote by $u_s$.  Here we have used $h$ to trivialize the density
bundle $| N^* \partial X|^{n-s}$.  We obtain
\begin{equation}\label{specialwit}
{\rm{pf} } \;
 \int_{ x >\epsilon}[|du_s|^2 -s(n-s)u_s^2]dv_g 
=-n\int_M S(s)1 dv_h.
\end{equation}
We consider  the limiting behaviour in this equation as $s\rightarrow n$.  
According to Theorem~\ref{t:2}, 
$\int_M S(s)1 dv_h \rightarrow c_{n/2}\int_M Qdv_h$.  We shall show that
the left hand side in (\ref{specialwit}) converges to $-nL/2$ as 
$s\rightarrow n$, thereby proving Theorem 3.

Since $u_s=\cP(s)1\rightarrow 1$ as $s\rightarrow n$, the integrand in the 
left hand side of (\ref{specialwit}) converges to 0 pointwise on $\Xint$.  
It follows that 
$$\int_{x>x_0}[|du_s|^2 -s(n-s)u_s^2]dv_g\rightarrow 0$$
for each fixed $x_0>0$.  So it suffices to consider
\begin{equation}\label{collar}
\int_{\ep<x<x_0}[|du_s|^2 -s(n-s)u_s^2]dv_g,
\end{equation}
for which we may use the
product identification \eqref{normalform}
coming from the representative metric $h$, provided
we choose $x_0$ small enough.  

By \eqref{eq:poiss'}, \eqref{eq:fexp}, \eqref{corr1}, and \eqref{Qs}, for
$s$ near $n$ we have  
\begin{equation}\label{uexp}
u_s = x^{n-s}(1+\sum_{k=1}^{n/2} c_{k,s}(n-s)Q_{k,s}x^{2k})
+x^sS(s)1 + \Oo(x^{n+3/4}),
\end{equation}
where the power $ n + 3/4 $ can be replaced by $ n + a $, $ a < 1 $, 
if we take $ s $ close to $ n$.  
In considering \eqref{uexp}, recall that 
$c_{k,s}$ is regular at $s=n$ for $1\leq k < n/2$, but that $c_{n/2,s}$ has 
a simple pole at $s=n$ with residue $c_{n/2}$.  Also, by Theorem 2,
$S(n)1=c_{n/2}Q=c_{n/2}Q_{n/2,n}$.  
Since $u_n=1$, the error term vanishes identically when
$s=n$ and in general is seen to be of the form $\Oo(|n-s|x^{n+3/4})$.
There is a similar bound when we differentiate the expansion.
Let $y^i$ denote local coordinates on $M$ and let $\pa$ denote any first
coordinate derivative $\pa_x$ or $\pa_{y^i}$; then we have 
\begin{equation}\label{duexp}
x\pa u_s=x\pa \left[ x^{n-s}(1+\sum_{k=1}^{n/2}
c_{k,s}(n-s)Q_{k,s}x^{2k}) +x^sS(s)1\right]+ \Oo(|n-s|x^{n+3/4}).
\end{equation}

We begin by considering the $u_s^2$ term in \eqref{collar}.  Upon squaring  
\eqref{uexp}, one obtains
$$
u_s^2=x^{2(n-s)}(1+\sum_{k=1}^{n/2} A_{k,s}x^{2k})
+2x^nS(s)1 + \Oo(x^{n+1/2}),
$$
where the coefficients $A_{k,s}$ are smooth functions on $M$, holomorphic
in $s$ near $s=n$, and satisfying $A_{k,n}=0$ for $1\leq k <n/2$ and 
$A_{n/2,n}=-2c_{n/2}Q$. 
Multiplying by $dv_g$ and using \eqref{volumeform} gives
$$
u_s^2dv_g=\left[ x^{n-2s-1}(1+\sum_{k=1}^{n/2} B_{k,s}x^{2k})
+2x^{-1}S(s)1 + \Oo(x^{-1/2})\right] dv_hdx,
$$
with coefficients $B_{k,s}$ again holomorphic in $s$, and with
$B_{n/2,n}=-2c_{n/2}Q + v^{(n)}$.
In order to evaluate $\lim_{s\rightarrow n} (n-s)\;{\rm{pf} } \,
\int_{\epsilon < x < x_0}u_s^2dv_g$, observe first that the $\Oo(x^{-1/2})$ 
error term is integrable, so its contribution vanishes upon letting 
$s\rightarrow n$.  Now
$$
{\rm{pf} } \;\int_{\epsilon < x < x_0} x^{n-2s+2k-1}B_{k,s}dv_hdx 
=x_0^{n-2s+2k}(n-2s+2k)^{-1}\int_M B_{k,s}dv_h.
$$
If $k< n/2$, then this approaches a finite limit as $s\rightarrow n$, 
so these terms also do not contribute.  The same is 
true  for the $x^{-1}S(s)1$
term.  Evaluating the limit for $k=n/2$ and recalling \eqref{Lform} then
gives 
\begin{equation}\label{pfu2}
\lim_{s\rightarrow n} (n-s)\;{\rm{pf} } \,
\int_{\epsilon < x < x_0}u_s^2dv_g = \frac12 \int_M B_{n/2,n}dv_h
=-c_{n/2}\int_M Q + L/2.
\end{equation}

For the derivative term, we have $|du_s|^2=(x\pa_x u_s)^2 +
h_x^{ij}(x\pa_{y^i}u_s)(x\pa_{y^j}u_s)$.  Consider first $(x\pa_x u_s)^2$.    
Expanding and squaring \eqref{duexp} gives
\begin{equation}\label{xdu2}
(x\pa_x u_s)^2=x^{2(n-s)}\sum_{k=0}^{n/2} A'_{k,s}x^{2k}
+2s(n-s)x^nS(s)1 + \Oo(|n-s|x^{n+1/2}),
\end{equation}
where $A'_{k,s}$ are smooth functions on $M$, holomorphic in $s$, 
satisfying 
\begin{equation}\label{coef}
\begin{split}
A'_{k,s}=\Oo(|n-s|^2), \;\; 0\leq k <n/2,\\
\lim_{s\rightarrow n}(n-s)^{-1}A'_{n/2,s}=-2nc_{n/2}Q.
\end{split}
\end{equation}
Multiplying \eqref{xdu2} by $dv_g$ gives
\[
\begin{split}
(x\pa_x u_s)^2 & dv_g=\\
&\left[ x^{n-2s-1}\sum_{k=0}^{n/2} B'_{k,s}x^{2k}
+2s(n-s)x^{-1}S(s)1 + \Oo(|n-s|x^{-1/2})\right] dv_h dx,
\end{split}
\]
where the $B'_{k,s}$ satisfy the same properties \eqref{coef} as the 
$A'_{k,s}$.
Integrating and evaluating the finite part and the limit as above yield
\begin{equation}\label{pfdx}
\lim_{s\rightarrow n} \;{\rm{pf} } \,
\int_{\epsilon < x < x_0}(x\pa_xu_s)^2dv_g.
=-nc_{n/2}\int_M Q\;. 
\end{equation}
Differentiation with respect to $y$ does not decrease the order of
vanishing in $x$, and because of this one finds by a similar calculation 
that 
\begin{equation}\label{pfdy}
\lim_{s\rightarrow n} \;{\rm{pf} } \,
\int_{\epsilon < x < x_0}h_x^{ij}(x\pa_{y^i}u_s)(x\pa_{y^j}u_s)dv_g
=0.
\end{equation}

Combining \eqref{pfu2}, \eqref{pfdx}, and \eqref{pfdy} gives
$$
\lim_{s\rightarrow n} \;{\rm{pf} } \,
\int_{\epsilon < x < x_0}[|du_s|^2-s(n-s)u_s^2] dv_g
= -nL/2
$$
as desired.
\stopthm

\end{document}